\documentstyle[12pt]{amsart}
\def\opn#1#2{\def#1{\operatorname{#2}}} % to make operators
\opn\chara{char}
\opn\length{\ell}
\opn\pd{pd}
\opn\rk{rk}
\opn\projdim{proj\,dim}
\opn\injdim{inj\,dim}
\opn\rank{rank}
\opn\depth{depth}
\opn\grade{grade}
\opn\height{height}
\opn\embdim{emb\,dim}
\opn\codim{codim}

\opn\Tr{Tr}
\opn\bigrank{big\,rank}
\opn\superheight{superheight}\opn\lcm{lcm}
\opn\trdeg{tr\,deg}%
\opn\reg{reg}
\opn\lreg{lreg}
\opn\skel{skel}
\opn\com{com}
\opn\div{div}
\opn\Div{Div}
\opn\cl{cl}
\opn\Cl{Cl}
\opn\Spec{Spec}
\opn\Supp{Supp}
\opn\supp{supp}
\opn\Sing{Sing}
\opn\Ass{Ass}
\opn\Ann{Ann}
\opn\Rad{Rad}
\opn\Soc{Soc}
\opn\Ker{Ker}
\opn\Coker{Coker}
\opn\Im{Im}
\opn\Hom{Hom}
\opn\Tor{Tor}
\opn\Ext{Ext}
\opn\End{End}
\opn\Aut{Aut}
\opn\id{id}

\opn\nat{nat}
\opn\pff{pf}%   \pf exists already
\opn\Pf{Pf}
\opn\GL{GL}
\opn\SL{SL}
\opn\mod{mod}
\opn\ord{ord}
\opn\aff{aff}
\opn\con{conv}
\opn\relint{relint}
\opn\st{st}
\opn\lk{lk}
\opn\cn{cn}
\opn\core{core}
\opn\vol{vol}
\opn\link{link}
\opn\star{star}
\opn\gr{gr}

\def\pot#1#2{#1[\kern-0.28ex[#2]\kern-0.28ex]}

\opn\dirlim{\underrightarrow{\lim}}
\opn\inivlim{\underleftarrow{\lim}}
\def\Implies{\ifmmode\Longrightarrow \else
     \unskip${}\Longrightarrow{}$\ignorespaces\fi}
\def\implies{\ifmmode\Rightarrow \else
     \unskip${}\Rightarrow{}$\ignorespaces\fi}
\def\iff{\ifmmode\Longleftrightarrow \else
     \unskip${}\Longleftrightarrow{}$\ignorespaces\fi}

\let\:=\colon
\newtheorem{Theorem}{Theorem}[section]
\newtheorem{Lemma}[Theorem]{Lemma}

\newtheorem{Question}[Theorem]{Question}
\let\epsilon\varepsilon
\let\phi=\varphi
\let\kappa=\varkappa
\def\qed{\ifhmode\textqed\fi
   \ifmmode\ifinner\quad\qedsymbol\else\dispqed\fi\fi}
\def\textqed{\unskip\nobreak\penalty50
    \hskip2em\hbox{}\nobreak\hfil\qedsymbol
    \parfillskip=0pt \finalhyphendemerits=0}
\def\dispqed{\rlap{\qquad\qedsymbol}}

\opn\initial{in}
\opn\inim{inm}
\opn\rev{rev}
\opn\Gin{Gin}
\opn\Lex{Lex}
\opn\Shift{Shift}
\opn\shift{shift}
\opn\rate{rate}
\opn\Mon{Mon}
\opn\lex{lex}
\opn\rev{rev}
\opn\red{red}
\opn\max{max}
\opn\min{min}
\opn\initial{in}
\opn\Ker{Ker}
\opn\GL{GL}
\opn\proj{proj}
\begin{document}
\title{$\Gin$ and $\Lex$ of certain monomial ideals}
\author{Satoshi Murai and Takayuki Hibi}
\date{}
\maketitle
\begin{abstract}
Let $A = K[x_1, \ldots, x_n]$ denote the polynomial ring 
in $n$ variables over a field $K$ of characteristic $0$
with each $\deg x_i = 1$.
Given arbitrary integers $i$ and $j$ with
$2 \leq i \leq n$ and $3 \leq j \leq n$,
we will construct a monomial ideal $I \subset A$ 
such that
(i) $\beta_k(I) < \beta_k(\Gin(I))$
for all $k < i$, 
(ii) $\beta_i(I) = \beta_i(\Gin(I))$,
(iii) $\beta_\ell(\Gin(I)) < \beta_\ell(\Lex(I))$
for all $\ell < j$ and
(iv) $\beta_j(\Gin(I)) = \beta_j(\Lex(I))$,
where $\Gin(I)$ is the
generic initial ideal 
of $I$ with respect to the reverse lexicographic order
induced by $x_1 > \cdots > x_n$
and where $\Lex(I)$ is the lexsegment ideal 
with the same Hilbert function as $I$. 
\end{abstract}

\section*{Introduction}
Let $A = K[x_1, \ldots, x_n]$ denote the polynomial ring 
in $n$ variables over a field $K$ of characteristic $0$
with each $\deg x_i = 1$.
Let $<_{\lex}$ (resp. $<_{\rev}$) denote the lexicographic  
(resp. reverse lexicographic) 
order on $A$ induced by 
the ordering $x_1 > \cdots > x_n$ of the variables.
Given a homogeneous ideal $I$ of $A$, we write $\Gin(I)$ for the
generic initial ideal (\cite[p. 348]{Eisenbud}) 
of $I$ with respect to $<_{\rev}$
and $\Lex(I)$ for the lexsegment ideal 
(\cite{Bigatti} and \cite{Hulett}) with the same Hilbert function 
as $I$. 
Let $\beta_i(I) = \dim_K \Tor^A_i(K,I)$ denote the $i$th 
Betti number of $I$ over $A$.
It is known that
\[
\beta_i(I) \leq \beta_i(\Gin(I)) \leq \beta_i(\Lex(I))
\]
for all $i \geq 0$. 
One has $\beta_i(I) = \beta_i(\Gin(I))$
for all $i$ if and only if $I$ is componentwise linear
(\cite{AramovaHerzogHibi}).  
One has $\beta_i(I) = \beta_i(\Lex(I))$
for all $i$ if and only if $I$ is Gotzmann (\cite{HerzogHibi}). 
On the other hand, the following facts are due to 
\cite[Corollary 2.7]{ConcaHerzogHibi}:

\begin{enumerate}
% \begin{itemize}
\item[(i)]
% \item
If $\beta_i(I) = \beta_i(\Gin(I))$ 
for some $i$,
then $\beta_k(I) = \beta_k(\Gin(I))$ 
for all $k \geq i$;
\item[(ii)]
% \item
If $\beta_i(I) = \beta_i(\Lex(I))$
for some $i$,
then $\beta_k(I) = \beta_k(\Lex(I))$
for all $k \geq i$.
% \end{itemize}
\end{enumerate}

\noindent
These behaviors of Betti numbers of $\Gin(I)$
and $\Lex(I)$ would naturally 
lead us to present the following 

\begin{Question}
\label{question}
Given arbitrary integers $1 \leq i \leq j \leq n$, 
does there exist a monomial ideal $I$ of
$A = K[x_1, \ldots, x_n]$ with the properties that
\begin{enumerate}
\item[(i)]
% (i) 
$\beta_k(I) < \beta_k(\Gin(I))$
for all $k < i$;
\item[(ii)]
%(ii) 
$\beta_i(I) = \beta_i(\Gin(I))$;
\item[(ii)]
% (iii) 
$\beta_\ell(I) < \beta_\ell(\Lex(I))$
for all $\ell < j$;
\item[(ii)]
% (iv) 
$\beta_j(I) = \beta_j(\Lex(I))$
{\em ?}
\end{enumerate}
\end{Question}

% It will turn out that the answer of Question \ref{question} 
% is affirmative.  
% In the case of $i < j$ a monomial ideal with the required
% properties will be constructed in Theorem \ref{Boston}.  
% On the other hand, in the case of $i = j$,
% the more detailed study will be achieved 
% in Theorem \ref{Sydney} together with Lemma \ref {lemma}
% and Example \ref{example}.

% To examine the monomial ideal arising in Theorem \ref{theorem} 
% will lead us to define the notion of 
% ``complementary lexsegment ideals.''
% In Theorem \ref{linear}
% we will discuss the problem when a complementary
% lexsegment ideal has a linear resolution.

The above question \ref{question} does ask the relation
between the Betti numbers of $I$ and $\Gin(I)$
together with the relation between 
the Betti numbers of $I$ and $\Lex(I)$.
In the present paper, however, the study of the relation
among the Betti numbers of $I$, $\Gin(I)$ and $\Lex(I)$
will be achieved.  Our goal is to show the following

\begin{Theorem}
\label{GinLex}
Let $A = K[x_1, \ldots, x_n]$ denote the polynomial ring 
in $n$ variables over a field $K$ of characteristic $0$
with each $\deg x_i = 1$.
Given arbitrary integers $i$ and $j$ with
$2 \leq i \leq n$ and $3 \leq j \leq n$,
there exists a monomial ideal $I \subset A$ 
such that
\begin{enumerate}
\item[(i)]
% (i) 
$\beta_k(I) < \beta_k(\Gin(I))$
for all $k < i$; 
\item[(ii)]
% (ii) 
$\beta_i(I) = \beta_i(\Gin(I))$;
\item[(iii)]
% (iii) 
$\beta_\ell(\Gin(I)) < \beta_\ell(\Lex(I))$
for all $\ell < j$;
\item[(iv)]
% (iv) 
$\beta_j(\Gin(I)) = \beta_j(\Lex(I))$.
\end{enumerate}
\end{Theorem}

A monomial ideal required in Theorem \ref{GinLex} 
will be given in Section $2$ 
in case of $2 \leq i < j \leq n$ and in Section $3$ 
in case of $3 \leq j \leq i \leq n$.
On the other hand, in Section $1$,
the reason why the assumption $3 \leq j$ is indispensable
in Theorem \ref{GinLex} will be explained.

\section{Betti numbers of strongly stable ideals}
Let $A = K[x_1, \ldots, x_n]$ denote the polynomial ring 
in $n$ variables over a field $K$ of characteristic $0$
with each $\deg x_i = 1$.
Recall that a monomial ideal $I$ of $A$ is {\em strongly stable}
if $u$ is a monomial belonging to $I$ and 
if $u$ is divided by $x_p$,
then $x_q u / x_p \in I$ for all $q < p$.
Since the base field $K$ is of characteristic $0$,
the generic initial ideal $\Gin(I)$ of an arbitrary homogeneous 
ideal of $A$ is strongly stable (\cite[Theorem 15.23]{Eisenbud}). 

Lemma \ref{lemma} below explains the reason  
why the assumption
$3 \leq j$ is indispensable in Theorem \ref{GinLex}.

\begin{Lemma}
\label{lemma}
Let $I \subset A$ be a strongly stable monomial ideal and
suppose that $\beta_i(I) = \beta_i(\Lex(I))$
for all $i \geq 2$.
Then $\beta_i(I) = \beta_i(\Lex(I))$ for all 
$i \geq 0$.
\end{Lemma}

\begin{pf}
Let $I$ be a strongly stable ideal.
We write $m_{\leq k}(I,j)$ for the number of
those monomials $u \in I$ of degree $j$ with
$m(u) \leq k$, where $m(u)$ is the largest integer $q$
for which $x_q$ divides $u$.
The computation of the Betti numbers of a strongly stable ideal
can be done by using the formura 
\cite[Proposition 2.3]{Bigatti} which says that
\begin{eqnarray}
\beta_{i,i+j}(I) = C - D,
\end{eqnarray}
where
\begin{eqnarray*}
C & = & m_{\leq n}(I,j) {n - 1 \choose i}, \\ 
D & = & \sum_{k=i}^{n-1} m_{\leq k}(I,j){k - 1 \choose i - 1}
+ \sum_{k=i+1}^{n} m_{\leq k}(I,j-1){k - 1 \choose i}.
\end{eqnarray*}
Since $\Lex(I)$ and $I$ have the same Hilbert function,
one has $m_{\leq n}(\Lex(I),j) = m_{\leq n}(I,j)$
for all $j$.  In addition, it is known
\cite[Theorem 2.1]{Bigatti} that 
$m_{\leq k}(\Lex(I),j) \leq m_{\leq k}(I,j)$
for all $j$ and $k$.
Thus, since $\beta_{i,i+j}(I) = \beta_{i,i+j}(\Lex(I))$
for all $i \geq 2$ and for all $j$, it follows that
$m_{\leq k}(\Lex(I),j) = m_{\leq k}(I,j)$
for all $j$ and for all $k \geq 2$.
On the other hand, since $I$ is strongly stable,
one has $m_{\leq 1}(I,j) = 1$ unless $I_j = (0)$,
where $I_j$ is the $j$th graded component of $I$. 
Hence $m_{\leq 1}(\Lex(I),j) = m_{\leq 1}(I,j)$
for all $j$.  Thus 
$m_{\leq k}(\Lex(I),j) \leq m_{\leq k}(I,j)$
for all $j$ and for all $k$.
\end{pf}

On the other hand, the reason why the assumption
$2 \leq i$ is indispensable in Theorem \ref{GinLex} 
is clear.  In fact, 
if $\beta_1(I) = \beta_1(\Gin(I))$, then
$\beta_i(I) = \beta_i(\Gin(I))$ for all $i \geq 0$.
To see why this is true, suppose that
$\beta_1(I) = \beta_1(\Gin(I))$.
Then
$\beta_i(I) = \beta_i(\Gin(I))$ for all $i \geq 1$.
Since $\beta_{i, i+j}(I) \leq \beta_{i, i +j}(\Gin(I))$
for all $i$ and for all $j$, one has
$\beta_{i, i+j}(I) = \beta_{i, i +j}(\Gin(I))$
for all $i \geq 1$ and for all $j$.
Since both $I$ and $\Gin(I)$ have the same Hilbert 
function, it follows that
$\beta_{0, j}(I) = \beta_{0, j}(\Gin(I))$
for all $j$.  Thus in particular
$\beta_0(I) = \beta_0(\Gin(I))$, as desired.

\section{Construction in the case of $2 \leq i < j \leq n$}
Theorem \ref{GinLex} is divided into 
Theorem \ref{Boston} and Theorem \ref{Sydney}.
A desired monomial ideal in the case of
$2 \leq i < j \leq n$   
will be given in Theorem \ref{Boston}
and that in the case of $3 \leq j \leq i \leq n$
will be given in Theorem \ref{Sydney}.

\begin{Theorem}
\label{Boston}
Let $A = K[x_1, \ldots, x_n]$ denote the polynomial ring 
in $n$ variables over a field $K$ of characteristic $0$
with each $\deg x_i = 1$.
Fix arbitrary integers $i$ and $j$ with 
$1 < i < j \leq n$ and $J$ the monomial ideal
of $A$ which is generated by those quadratic monomials 
$x_p x_q$, $1 \leq p \leq q \leq n$, with
$x_{i-1} x_{j} <_{\lex} x_p x_q$.
Suppose that $I$ is the monomial ideal
\[
I = (x_1, \ldots, x_n)^3 + J + (x_n^2) 
\]
of $A$.  Then one has 
\begin{enumerate}
\item[(i)]
$\beta_k(I) < \beta_k(\Gin(I))$
for all $k < i$;
\item[(ii)]
$\beta_i(I) = \beta_i(\Gin(I))$;
\item[(iii)]
$\beta_\ell(\Gin(I)) < \beta_\ell(\Lex(I))$
for all $\ell < j$;
\item[(iv)]
$\beta_j(\Gin(I)) = \beta_j(\Lex(I))$.
\end{enumerate}
\end{Theorem}

\begin{pf}
(\,{\em First Step}\,)
% First of all, in general, 
Given a monomial ideal $L$ of 
$A = K[x_1, \ldots, x_n]$ % which is 
generated by
quadratic monomials, we introduce  
a finite graph $\Gamma(L)$ on the vertex set
\[
V = \{ 1, \ldots, n, 1', \ldots, n' \}
\]
whose edge set $E(\Gamma(L))$ consists of 
those edges
\begin{itemize}
\item [(i)]
$\{ p, q \}$ with $1 \leq p < q \leq n$
such that $x_p x_q \not\in L$;
\item [(ii)]
$\{ p, q' \}$ with $1 \leq p \leq n$,
$1 \leq q \leq n$ and $p \neq q$;
\item [(iii)]
$\{ p, p' \}$ with $1 \leq p \leq n$
such that $x_p^2 \not\in L$.
\end{itemize}
If $W \subset V$, then we write $\Gamma(L)_W$
for the induced subgraph of $\Gamma(L)$ on $W$.
Let $\delta(\Gamma(L)_W)$ denote the number of connected
component of $\Gamma(L)_W$.

By using Hochster's formula
\cite[Theorem 5.1.1]{BrunsHerzog}
together with the polarization technique
\cite[Lemma 4.2.16]{BrunsHerzog},
it follows that the Betti number
$\beta_{k,k+2}(L) = \dim_K [\Tor^A_k(K,I)]_{k+2}$
can be computed by the formula
\begin{eqnarray}
\beta_{k,k+2}(L) 
= \sum_{W \subset V, \, \, |W| = k + 2} 
(\delta(\Gamma(L)_W) - 1).
\end{eqnarray}

\medskip

\noindent
(\,{\em Second Step}\,)
Let $\Gin_{<_{\lex}}(I)$ denote the generic initial ideal 
of $I$ with respect to 
the lexicographic order $<_{\lex}$.
We claim
\begin{eqnarray*}
\Gin(I) & = & (x_1, \ldots, x_n)^3 + J + (x_i^2); \\
\Lex(I) & = & \Gin_{<_{\lex}}(I)
\, = \, (x_1, \ldots, x_n)^3 + J + (x_{i-1} x_{j}).
\end{eqnarray*}

Assume that the general linear group
$\GL(n;K)$ acts linearly on $A$.
Let $\psi \in \GL(n;K)$.  Then
$\Gin(\psi(I)) = \Gin(I)$ 
and $\Gin_{<_{\lex}}(\psi(I)) = \Gin_{<_{\lex}}(I)$.

\begin{itemize}
\item
Let $\varphi \in \GL(n;K)$ be defined by
$\varphi(x_i) = x_i$ for all $i \neq n$
and $\varphi(x_n) = x_i + x_n$.
Then the initial ideal 
$\initial_{<_{\rev}}(\varphi(I))$
of $\varphi(I)$ with respect to $<_{\rev}$ 
coincides with the strongly stable ideal
$(x_1, \ldots, x_n)^3 + J + (x_i^2)$.
Hence $\Gin(\initial_{<_{\rev}}(\varphi(I)))
= (x_1, \ldots, x_n)^3 + J + (x_i^2)$.
Since $\initial_{<_{\rev}}(\varphi(I))$ 
contains all monomials 
$x_p x_q$, $1 \leq p \leq q \leq n$,
with $x_i^2 \leq_{\rev} x_p x_q$,
it follows from \cite[Corollary 1.6]{Conca} that
$x_i^2 \in \Gin(\varphi(I))$.
Moreover, since $J$ is strongly stable and
$J \subset \varphi(I)$, one has
$J = \Gin(J) \subset \Gin(\varphi(I))$.
Thus $(x_1, \ldots, x_n)^3 + J + (x_i^2)
\subset \Gin(\varphi(I)) = \Gin(I)$.
Hence $\Gin(I) = (x_1, \ldots, x_n)^3 + J + (x_i^2)$,
as desired.

\smallskip

\item
Let $\varphi \in \GL(n;K)$ be defined by
$\varphi(x_i) = x_i$ for all $i \neq n$
and $\varphi(x_n) = x_{i-1} + x_{j} + x_n$.
Then the initial ideal 
$\initial_{<_{\lex}}(\varphi(I))$
of $\varphi(I)$ with respect to $<_{\lex}$ 
coincides with the strongly stable ideal
$(x_1, \ldots, x_n)^3 + J + (x_{i-1} x_{j})$.
Hence $\Gin_{<_{\lex}}(\initial_{<_{\lex}}(\varphi(I)))
= (x_1, \ldots, x_n)^3 + J + (x_{i-1} x_{j})$.
Since $\initial_{<_{\lex}}(\varphi(I))$ 
contains all monomials 
$x_p x_q$, $1 \leq p \leq q \leq n$,
with $x_{i-1} x_{j} \leq_{\lex} x_p x_q$,
it follows from \cite[Corollary 1.6]{Conca} that
$x_{i-1} x_{j} \in \Gin_{<_{\lex}}(\varphi(I))$.
Moreover, since $J$ is strongly stable and
$J \subset \varphi(I)$, one has
$J = \Gin_{<_{\lex}}(J) \subset 
\Gin_{<_{\lex}}(\varphi(I))$.
Thus $(x_1, \ldots, x_n)^3 + J + (x_{i-1} x_{j})
\subset \Gin_{<_{\lex}}(\varphi(I)) 
= \Gin_{<_{\lex}}(I)$.
Hence $\Gin_{<_{\lex}}(I) 
= (x_1, \ldots, x_n)^3 + J + (x_{i-1} x_{j})$.
Since the ideal $J + (x_{i-1} x_{j})$
is lexsegment, one has
$\Lex(I) = (x_1, \ldots, x_n)^3 + J + (x_{i-1} x_{j})$.
Hence 
$\Lex(I) = \Gin_{<_{\lex}}(I) 
= (x_1, \ldots, x_n)^3 + J + (x_{i-1} x_{j})$,
as desired.
\end{itemize}

\medskip

\noindent
(\,{\em Third Step}\,)
We compute
\begin{eqnarray*}
\beta_{k,k+2}(I) & = & \beta_{k,k+2}(J + (x_n^2)); \\ 
\beta_{k,k+2}(\Gin(I)) & = & \beta_{k,k+2}(J + (x_i^2)); \\
\beta_{k,k+2}(\Lex(I)) 
& = & \beta_{k,k+2}(J + (x_{i-1} x_{j})),
\end{eqnarray*}
based on the formula $(1)$
together with the combinatorics on the finite graphs
$\Gamma(J),
\Gamma(J + (x_n^2)),
\Gamma(J + (x_i^2))$ and
$\Gamma(J + (x_{i-1} x_{j}))$
with
\begin{eqnarray*}
E(\Gamma(J + (x_n^2))) 
& = & E(\Gamma(J)) \setminus \{ \{ n, n' \} \}; \\
E(\Gamma(J + (x_i^2))) 
& = & E(\Gamma(J)) \setminus \{ \{ i, i' \} \}; \\
E(\Gamma(J + (x_{i-1} x_{j}))) 
& = & E(\Gamma(J)) \setminus \{ \{ i-1, j \} \}.
\end{eqnarray*}

\begin{itemize}
\item
Let $W \subset V$ with $|W| = k + 2$ 
and suppose that
$\delta(\Gamma(J + (x_n^2))_W) >
\delta(\Gamma(J)_W)$.
Then
(i) both $n$ and $n'$ belong to $W$,
(ii) no connected component of 
$\Gamma(J + (x_n^2))_W$ contains both $n$ and $n'$,
and (iii)
$\delta(\Gamma(J + (x_n^2))_W) =
\delta(\Gamma(J)_W) + 1$.
Since $\{ \alpha , n' \}$ is an edge of 
$\Gamma(J)$ for all $\alpha \in V$
with $\alpha \neq n'$
and since $\{ \beta , n \}$ is an edge of 
$\Gamma(J)$ if and only if
$\beta \in V \setminus \{ 1, 2, \ldots, i - 2 \}$,
it follows that 
$W \setminus \{ n, n' \} \subset \{ 1, 2, \ldots, i - 2 \}$.
Hence the number of subsets $W \subset V$ with
$|W| = k + 2$ 
such that 
$\delta(\Gamma(J + (x_n^2))_W) >
\delta(\Gamma(J)_W)$
is ${ i - 2 \choose k }$.

\smallskip

\item
Let $W \subset V$ with $|W| = k + 2$ 
and suppose that
$\delta(\Gamma(J + (x_i^2))_W) >
\delta(\Gamma(J)_W)$.
Then
(i) both $i$ and $i'$ belong to $W$,
(ii) no connected component of 
$\Gamma(J + (x_n^2))_W$ contains both $i$ and $i'$,
and (iii)
$\delta(\Gamma(J + (x_i^2))_W) =
\delta(\Gamma(J)_W) + 1$.
Since $\{ \alpha , i' \}$ is an edge of 
$\Gamma(J)$ for all $\alpha \in V$
with $\alpha \neq i'$
and since $\{ \beta , i \}$ is an edge of 
$\Gamma(J)$ if and only if
$\beta \in V \setminus \{ 1, 2, \ldots, i - 1 \}$,
it follows that 
$W \setminus \{ i, i' \} \subset \{ 1, 2, \ldots, i - 1 \}$.
Hence the number of subsets $W \subset V$ with
$|W| = k + 2$ 
such that 
$\delta(\Gamma(J + (x_i^2))_W) >
\delta(\Gamma(J)_W)$
is ${ i - 1 \choose k }$.

\smallskip

\item
Let $W \subset V$ with $|W| = k + 2$ 
and suppose that
$\delta(\Gamma(J + (x_{i-1} x_{j}))_W) >
\delta(\Gamma(J)_W)$.
Then
(i) both $i - 1$ and $j$ belong to $W$,
(ii) no connected component of 
$\Gamma(J + (x_{i-1} x_{j}))_W$ contains both 
$i - 1$ and $j$,
and (iii)
$\delta(\Gamma(J + (x_{i-1} x_{j}))_W) =
\delta(\Gamma(J)_W) + 1$.
Since $\{ \alpha , i - 1 \}$ is an edge of 
$\Gamma(J)$ if and only if 
\[
\alpha \in 
\{ j, j + 1, \ldots, n \}
\bigcup
(\{ 1', 2', \ldots, n' \}
\setminus \{ (i - 1)' \})
\]
and since $\{ \beta , j \}$ is an edge of 
$\Gamma(J)$ if and only if
\[
\beta \in (\{ i - 1, i, \ldots, n  \} \setminus \{ j - 1 \})
\bigcup \{ 1', 2', \ldots, n' \},
\]
it follows that 
\[
W \setminus \{ i, i' \} \subset \{ 1, 2, \ldots, i - 2,
(i - 1)', i, i + 1, \ldots, j - 1 \}.
\]
Hence the number of subsets $W \subset V$ with
$|W| = k + 2$ 
such that
$\delta(\Gamma(J + (x_{i-1} x_{j}))_W) >
\delta(\Gamma(J)_W)$
is ${ j - 1 \choose k }$.
\end{itemize}

\noindent
Hence 
\begin{eqnarray*}
\beta_{k,k+2}(I)
& = & \beta_{k,k+2}(J) + { i - 2 \choose k }; \\
\beta_{k,k+2}(\Gin(I))
& = & \beta_{k,k+2}(J) + { i - 1 \choose k }; \\
\beta_{k,k+2}(\Lex(I))
& = & \beta_{k,k+2}(J) + { j - 1 \choose k }. 
\end{eqnarray*}
Thus in particular
\begin{enumerate}
\item[(i)]
$\beta_{k,k+2}(I) < \beta_{k,k+2}(\Gin(I))$
for all $k < i$;
\item[(ii)]
$\beta_{k,k+2}(I) = \beta_{k,k+2}(\Gin(I))$
for all $k \geq i$;
\item[(iii)]
$\beta_{\ell,\ell+2}(\Gin(I)) 
< \beta_{\ell,\ell+2}(\Lex(I))$
for all $\ell < j$;
\item[(iv)]
$\beta_{\ell,\ell+2}(\Gin(I)) 
= \beta_{\ell,\ell+2}(\Lex(I))$
for all $\ell \geq j$.
\end{enumerate}

\medskip

\noindent
(\,{\em Fourth Step}\,)
Since the regularity of each of 
$I$, $\Gin(I)$ and $\Lex(I)$ is $3$, it follows that 
\[
\beta_{p,p+q}(I) = \beta_{p,p+q}(\Gin(I)) 
= \beta_{p,p+q}(\Lex(I)) = 0
\]
for all $q > 3$.     
The cancellation principle \cite[Corollary 1.21]{Green}
now guarantees that
\[
\beta_{k,k+2}(\Gin(I)) - \beta_{k,k+2}(I)
= 
\beta_{(k-1),(k-1)+3}(\Gin(I)) - \beta_{(k-1),(k-1)+3}(I)
\]
for all $k$.  Since 
$\beta_{k,k+2}(I) = \beta_{k,k+2}(\Gin(I))$
for all $k \geq i$, 
one has
\[
\beta_{(k-1),(k-1)+3}(I) 
= \beta_{(k-1),(k-1)+3}(\Gin(I))
\]
for all $k - 1 \geq i - 1$.
Hence
$\beta_{k}(I) = \beta_{k}(\Gin(I))$
for all $k \geq i$.
Similarly, since $\Lex(I) = \Gin_{<_{\lex}}(I)$,
by using the cancellation principle 
to $\Gin_{<_{\lex}}(I)$,
one has   
$\beta_{\ell}(I) = \beta_{\ell}(\Lex(I))$
for all $\ell \geq j$, as reqired.
\end{pf}

\section{Construction in the case of $3 \leq j \leq i \leq n$}
A monomial ideal
in the case of $3 \leq j \leq i \leq n$
in Theorem \ref{GinLex}
will be given in Theorem \ref{Sydney}.
Note that in Theorem \ref{Sydney} we use
$i + 1$ and $j + 1$ instead of $i$ and $j$, so that 
we work with fixing arbitrary integers $i$ and $j$ 
with $2 \leq j \leq i < n$. 
  
\begin{Theorem}
\label{Sydney}
Let $A = K[x_1, \ldots, x_n]$ denote the polynomial ring 
in $n$ variables over a field $K$ of characteristic $0$
with each $\deg x_i = 1$.
Fix arbitrary integers $i$ and $j$ 
with $2 \leq j \leq i < n$.
Let $H$ denote the monomial ideal
of $A$ which is generated by those quadratic monomials 
$x_p x_q$, $1 \leq p \leq q \leq n$, with
$x_{j-1} x_{j} \leq_{\lex} x_p x_q$.
Let $G$ denote the monomial ideal
of $A$ which is generated by those quadratic monomials 
$x_p x_q$, $2 \leq p \leq q \leq n$, 
with
$x_{i}^2 \leq_{\lex} x_p x_q$.
Suppose that $I$ is the monomial ideal
\[
I = x_1(H + (x_n^2)) + x_1(x_1, \ldots, x_n)^3 
+ x_2^2(G + (x_n^2)) + (x_1, \ldots, x_n)^5 
\]
of $A$.  Then one has 
\begin{enumerate}
\item[(i)]
% (i) 
$\beta_k(I) < \beta_k(\Gin(I))$
for all $k \leq i$; 
\item[(ii)]
% (ii) 
$\beta_{i+1}(I) = \beta_{i+1}(\Gin(I))$;
\item[(iii)]
% (iii) 
$\beta_\ell(\Gin(I)) < \beta_\ell(\Lex(I))$
for all $\ell \leq j$;
\item[(iv)]
% (iv) 
$\beta_{j+1}(\Gin(I)) = \beta_{j+1}(\Lex(I))$.
\end{enumerate}
\end{Theorem}

\begin{pf}
(\,{\em First Step}\,)
\,
First, we discuss the relation between $\beta_k(\Gin(I))$ 
and $\beta_k(\Lex(I))$.
The same techniques as in Second Step of
the proof of Theorem \ref{Boston}
yields that 
\begin{eqnarray*}
\Gin(I) & = & 
x_1(H + (x_j^2)) + x_1(x_1, \ldots, x_n)^3 
+ x_2^2(G + (x_{i}x_{i+1})) + (x_1, \ldots, x_n)^5;
\\
\Lex(I) & = & 
x_1(H + (x_{j-1}x_{j+1})) + x_1(x_1, \ldots, x_n)^3 
+ x_2^2(G + (x_ix_{i+1})) + (x_1, \ldots, x_n)^5.
\end{eqnarray*}
Since
$\Gin(I)_d = \Lex(I)_d$ for all $d \geq 4$,
one has
$m_{\leq k}(\Gin(I),d) = m_{\leq k}(\Lex(I),d)$
for all $k$ and for all $d \geq 4$.
In addition, 
$m_{\leq k}(\Gin(I),3) = m_{\leq k}(\Lex(I),3)$
for all $k \neq j$ and
$m_{\leq j}(\Gin(I),3) < m_{\leq j}(\Lex(I),3)$.
It then follows from the formula $(2)$ that
\begin{eqnarray*}
\beta_{j+1}(\Gin(I)) 
& = & \beta_{j+1}(\Lex(I));
\\
\beta_{k}(\Gin(I)) 
& < & \beta_{k}(\Lex(I))
\, \, \, \, \, \, \, \, \, \, 
\text{for all} \, \, \, \, \,  
k \leq j.
\end{eqnarray*}

\medskip

\noindent
(\,{\em Second Step}\,)
\,
Let ${\tilde I} = x_1 H$.
We claim that
\begin{eqnarray*}
\beta_{k, k+3}(I) 
& = & \beta_{k, k+3}({\tilde I}) + {j - 2 \choose k};
\\
\beta_{k, k+3}(\Gin(I)) 
& = & \beta_{k, k+3}({\tilde I}) + {j - 1 \choose k}.
\end{eqnarray*}

Let $I_{\leq d}$ denote the ideal generated by
those monomials $u \in I$ with $\deg u \leq d$.
Recall that 
$\beta_{i,i+d}(I) = \beta_{i,i+d}(I_{\leq d})$
(\cite[Lemma 1.2]{HerzogHibi}).
Thus $\beta_{k, k+3}(I) = \beta_{k,k+3}(I_{\leq 3})$ and 
$\beta_{k, k+3}(\Gin(I)) = \beta_{k,k+3}(\Gin(I)_{\leq 3})$.
Since
$I_{\leq 3} = x_1(H + (x_n^2))$ and
$\Gin(I)_{\leq 3} = x_1(H + (x_j^2))$,
one has
$\beta_{k,k+3}(I_{\leq 3})
= \beta_{k,k+2}((H + (x_n^2))$ and
$\beta_{k,k+3}(\Gin(I)_{\leq 3})
= \beta_{k,k+2}((H + (x_j^2))$.
It follows from the same computation as in
Third Step of the proof of Thorem \ref{Boston} 
that
\begin{eqnarray*}
\beta_{k, k+3}(I) 
= \beta_{k,k+2}((H + (x_n^2))
= \beta_{k,k+2}(H) + {j - 2 \choose k}
= \beta_{k,k+3}({\tilde I}) + {j - 2 \choose k};
\\
\beta_{k, k+3}(\Gin(I)) 
= \beta_{k,k+2}((H + (x_j^2))
= \beta_{k,k+2}(H) + {j - 1 \choose k}
= \beta_{k,k+3}({\tilde I}) + {j - 1 \choose k}.
\end{eqnarray*}

\medskip

\noindent
(\,{\em Third Step}\,)
\,
We now turn to the computation of
$\beta_{k, k+4}(I)$ and $\beta_{k, k+4}(\Gin(I))$.
Let 
\[
J = x_1(H + (x_j^2)) + x_1(x_1, \ldots, x_n)^3 + x_2^2 G;
\]
\[
{\tilde J} = x_1(H + (x_n^2)) + x_1(x_1, \ldots, x_n)^3 + x_2^2 G.
\]
We claim that
\begin{eqnarray}
\beta_{k, k+4}(I) 
& = & \beta_{k, k+4}(J) + {i - 1 \choose k}
- {j - 1 \choose k + 1} + {j - 2 \choose k + 1};
\\
\beta_{k, k+4}(\Gin(I)) 
& = & \beta_{k, k+4}(J) + {i \choose k}.
\end{eqnarray}

\medskip

\noindent
{\bf (\,3\,.\,1\,)} \,
Since 
$J_{\leq 3} = \Gin(I)_{\leq 3}$ and
${\tilde J}_{\leq 3} = I_{\leq 3}$,
one has
$\beta_{k,k+3}(J) = \beta_{k,k+3}(\Gin(I))$
and
$\beta_{k,k+3}(\tilde J) = \beta_{k,k+3}(I)$.
Let $J_{\geq d}$ denote the ideal generated by those
monomials $u \in J$ with $\deg u \geq d$. 
Since $J_{\geq 4}$ and ${\tilde J}_{\geq 4}$ 
are strongly stable, the regularity of 
each of $J$ and ${\tilde J}$ is $4$.
Thus, since $\Gin({\tilde J}) = J$, 
the cancellation principle 
\cite[Corollary 1.21]{Green}
yields that
\begin{eqnarray}
\beta_{k,k+4}(J) - \beta_{k,k+4}({\tilde J})
= \beta_{k+1,k+1+3}(J) - \beta_{k+1,k+1+3}({\tilde J}).
\end{eqnarray}
By virtue of Second Step the right-hand side of $(5)$
is equal to ${j - 1 \choose k + 1} - {j - 2 \choose k + 1}$.
Thus
\begin{eqnarray}
\beta_{k,k+4}(J) = \beta_{k,k+4}({\tilde J})
+ {j - 1 \choose k + 1} - {j - 2 \choose k + 1}.
\end{eqnarray}

\medskip

\noindent
{\bf (\,3\,.\,2\,)} \,  
We now show the equality
\begin{eqnarray}
\beta_{k,k+4}(I) = \beta_{k,k+4}({\tilde J})
+ {i - 1 \choose k}.
\end{eqnarray}
Since $I_{\leq 4} = {\tilde J} + (x_2^2 x_n^2)$,
the shot exact sequence
\[
0 \longrightarrow
{\tilde J}
\longrightarrow
I_{\leq 4} = {\tilde J} + (x_2^2 x_n^2)
\longrightarrow
({\tilde J} + (x_2^2 x_n^2)) / {\tilde J}
\longrightarrow 0
\]
yields the long exact sequence
\begin{eqnarray*}
& & \hspace{6cm} 
\cdots \longrightarrow
[\Tor_{k+1}^A(K,({\tilde J} + (x_2^2 x_n^2)) / {\tilde J})]_{k+1+3}
\\
& \longrightarrow &
[\Tor_{k}^A(K,{\tilde J})]_{k+4}
\longrightarrow
[\Tor_{k}^A(K,I_{\leq 4})]_{k+4}
\longrightarrow
[\Tor_{k}^A(K,({\tilde J} + (x_2^2 x_n^2)) / {\tilde J})]_{k+4}
\\
& \longrightarrow &
[\Tor_{k-1}^A(K,{\tilde J})]_{k-1+5}
\longrightarrow \cdots.
\end{eqnarray*}
Since $(({\tilde J} + (x_2^2 x_n^2)) / {\tilde J})_3 = 0$, 
one has
$[\Tor_{k+1}^A(K,({\tilde J} + (x_2^2 x_n^2)) / {\tilde J})]_{k+1+3} 
= 0$.
Since the regularity of ${\tilde J}$ is $4$, one has
$[\Tor_{k-1}^A(K,{\tilde J})]_{k-1+5} = 0$.
Thus the above long exact sequence turns out to be
\begin{eqnarray*}
0 \longrightarrow 
[\Tor_{k}^A(K,{\tilde J})]_{k+4}
\longrightarrow
[\Tor_{k}^A(K,I_{\leq 4})]_{k+4}
\longrightarrow
[\Tor_{k}^A(K,({\tilde J} + (x_2^2 x_n^2)) / {\tilde J})]_{k+4}
\longrightarrow 0.
\end{eqnarray*}
In particular
\[
\beta_{k,k+4}(I_{\leq 4})) 
= \beta_{k,k+4}({\tilde J}) 
+ \beta_{k,k+4}(({\tilde J} + (x_2^2 x_n^2)) / {\tilde J}).
\]
Since $\beta_{k,k+4}(I) = \beta_{k,k+4}(I_{\leq 4})$,
to show the equality $(7)$, what we must prove is
\begin{eqnarray}
\beta_{k,k+4}(({\tilde J} + (x_2^2 x_n^2)) / {\tilde J})
= {i - 1 \choose k}.
\end{eqnarray}
Let ${\tilde G}$ denote the ideal generated by those
quadratic monomials $x_p x_q$, $1 \leq p < q \leq n$,
with $x_i^2 \leq_{\lex} x_px_q$.
A routine computation shows that
\[
({\tilde J} + (x_2^2 x_n^2)) / {\tilde J}
\cong x_2^2(({\tilde G} + (x_n^2))/{\tilde G}).
\]
Thus in particular
\begin{eqnarray}
\beta_{k,k+4}(({\tilde J} + (x_2^2 x_n^2)) / {\tilde J})
= \beta_{k,k+2}(({\tilde G} + (x_n^2))/{\tilde G}).
\end{eqnarray}
Again, the short exact sequence
\[
0 \longrightarrow
{\tilde G}
\longrightarrow
{\tilde G} + (x_n^2)
\longrightarrow
({\tilde G} + (x_n^2)) / {\tilde G}
\longrightarrow 0
\]
yields the long exact sequence
\begin{eqnarray*}
& & \hspace{6cm} 
\cdots \longrightarrow
[\Tor_{k+1}^A(K,({\tilde G} + (x_n^2)) / {\tilde G})]_{k+1+1}
\\
& \longrightarrow &
[\Tor_{k}^A(K,{\tilde G})]_{k+2}
\longrightarrow
[\Tor_{k}^A(K,{\tilde G} + (x_n^2))]_{k+2}
\longrightarrow
[\Tor_{k}^A(K,({\tilde G} + (x_n^2)) / {\tilde G})]_{k+2}
\\
& \longrightarrow &
[\Tor_{k-1}^A(K,{\tilde G})]_{k-1+3}
\longrightarrow \cdots.
\end{eqnarray*}
Since $({\tilde G} + (x_n^2)) / {\tilde G})_1 = 0$, 
one has
$[\Tor_{k+1}^A(K,({\tilde G} + (x_n^2)) / {\tilde G})]_{k+1+1}
= 0$.
Since the regularity of ${\tilde G}$ is $2$, one has
$[\Tor_{k-1}^A(K,{\tilde G})]_{k-1+3} = 0$.
Thus
\begin{eqnarray}
\beta_{k,k+2}(({\tilde G} + (x_n^2)) / {\tilde G})
= 
\beta_{k,k+2}({\tilde G} + (x_n^2)) 
- \beta_{k,k+2}({\tilde G}).
\end{eqnarray}
Again, the same computation as in Third Step 
of the proof of Theorem \ref{Boston} says that
\begin{eqnarray}
\beta_{k,k+2}({\tilde G} + (x_n^2)) 
- \beta_{k,k+2}({\tilde G})
= {i - 1 \choose k}.
\end{eqnarray}
The equalities $(9)$ and $(10)$ together with $(11)$ now 
yield the desired equality $(8)$. 

\medskip

\noindent
{\bf (\,3\,.\,3\,)}
\,
The first equality $(3)$ in our claim follows from
the equalities $(6)$ and $(7)$.
On the other hand, since $\Gin(I)$ and $J$ are strongly stable,
the second equality $(4)$ in our claim follows from 
the formula \cite[Corollary 3.6 (a)]{Herzog}
obtained by Eliahou and Kervaire \cite{EliahouKervaire}.

\medskip

\noindent
(\,{\em Fourth Step}\,)
\,
By virtue of Second Step and Third Step, 
it follows that
\begin{eqnarray*}
\beta_{k, k+3}(I) & = & \beta_{k, k+3}(\Gin(I))
\, \, \, \, \, \, \, \, \, \, 
\text{for all} \, \, \, \, \,  
k \geq j;
\\ 
\beta_{k, k+4}(I) & = & \beta_{k, k+4}(\Gin(I))
\, \, \, \, \, \, \, \, \, \,
\text{for all} \, \, \, \, \,  
k \geq i + 1.
\end{eqnarray*}
The cancellation principle \cite[Corollary 1.21]{Green}
then guarantees that
\[
\beta_{k, k+5}(I) = \beta_{k, k+5}(\Gin(I))
\, \, \, \, \, \, \, \, \, \, 
\text{for all} \, \, \, \, \,  
k \geq i.
\]
Since the regularity of each of $I$ and $\Gin(I)$ is $5$,
one has
$\beta_k(I) = \beta_k(\Gin(I))$
for all $k \geq i + 1$.
Again, by virtue of Third Step, it follows that
\[
\beta_{i,i+4}(I) = \beta_{i,i+4}(J) < \beta_{i,i+4}(J) + 1
= \beta_{i,i+4}(\Gin(I)).
\]
In particular $\beta_i(I) < \beta_i(\Gin(I))$.
Hence $\beta_k(I) < \beta_k(\Gin(I))$ for all $k \leq i$.
\end{pf}

\bigskip

{\small
\noindent
Satoshi Murai \\
Department of Pure and Applied Mathematics\\
Graduate School of Information Science and Technology\\
Osaka University \\
Toyonaka, Osaka 560-0043, Japan\\
E-mail:s-murai@@ist.osaka-u.ac.jp

\bigskip

\noindent
Takayuki Hibi\\
Department of Pure and Applied Mathematics\\
Graduate School of Information Science and Technology\\
Osaka University \\
Toyonaka, Osaka 560-0043, Japan\\
E-mail:hibi@@math.sci.osaka-u.ac.jp

}

\end{document}